\documentclass[A4paper,12pt]{article}
\usepackage{latexsym}
\usepackage{mathrsfs}
\usepackage{amssymb}
\usepackage{amscd}
\usepackage[dvips]{graphicx}  
\usepackage{xcolor}
\usepackage{color}

\newtheorem{definition}{Definition}
\newtheorem{theorem}{Theorem}

\newtheorem{lemma}{Lemma}

\newenvironment{proof}[1][Proof]{\textbf{#1.} }{\ \rule{0.5em}{0.5em}}
\date{}

\long\def\symbolfootnote[#1]#2{\begingroup%
	\def\thefootnote{$\;$}\footnote[#1]{$^*$#2}\endgroup}
\begin{document}
	
	\title{Nonmeasurable sets in tree structures \\and Ellentuck topology}
	\author{Joanna Jureczko}
	
		\maketitle
	
	\symbolfootnote[2]{Mathematics Subject Classification (2010): 54C30, 03E05, 03E40, 28A20.  
		
		\hspace{0.2cm}
		}
	
	\begin{abstract} 
		The aim of this paper is to provide the results that answer the Kuratowski problem posed in 1935 concerning the existence of nonmeasurable sets. The Kuratowski problem was considered for partitions, here we provide a generalization to point-finite families for tree structures and structures in Ellentuck space. The main result of this paper generalizes previous results in the topic of nonmeasurable sets.  
	\end{abstract}
	
	\section{Introduction}
	
	The problem of the existence of nonmeasurable sets for various structures is a topic that has been considered for a long time. Examples of such sets are Bernstein's or Vitali's. Over the years, there have been attempts to construct or demonstrate that such sets exist. The growing number of papers on the subject of the existence of nonmeasurable sets proves that this topic is still alive. However, so far only a small number of significantly equal nonmeasurable sets have been published. One of such sets is the one proposed by Kuratowski in 1935, and in fact the beginning of research in this direction has his question posed in the paper \cite{KK1}.
	
	In 1935  K. Kuratowski  posed the problem (\cite{KK1}) which is equivalent to non-existence of Kuratowski partitions in the term of functions. It was found to be equivalent to the following result, (see \cite{KK2}). 
	\\ 
	\\
	\textbf{Theorem (Kuratowski)}
	\textit{(CH) If $\mathcal{F}$  is a partition of the interval $[0,1]$ of cardinality at most $\omega_1$ into meager sets then there exists a family $\mathcal{A} \subseteq \mathcal{F}$ such that $\bigcup \mathcal{A}$ has not the Baire property.}
	\\
	
	As there is a duality between the concepts of "having the Baire property" and "being nonmeasurable", through the Kuratowski theorem we obtain the existence of an nonmeasurable set.
	
	The problem of the existence of nonmeasurable sets has been intensively researched since the 1970s' of the last century with using different methods of proving. For example R. Solovay in \cite{RS} and L. Bukovsky in \cite{LB} proved this result without (CH) assumption using the forcing methods and natural embeddings in their proofs.
	In \cite{EFK}, the authors gave a generalization of this result for a class of spaces of weight $\leq 2^\omega$. 
	Moreover this problem was also considered in \cite{BCGR} in special case and in \cite{LS}.
	
	In  \cite{BCGR} the problem of the existence of nonmeasurable sets has been considered for so called point-finite families instead of partitions.
	More precisely,  the authors showed the following result
		\\\\
		\textbf{Theorem(\cite{BCGR})}
		\textit{	Let $\mathbb{R}$ be the set of reals and let $\mathcal{B}$ be a $\sigma$-algebra of Borel sets in $\mathbb{R}$. Let $I$ be a $\sigma$-ideal on $\mathbb{R}$ such that $I \cap \mathcal{B}$ is a base for $I$. If $\{A_s \colon s \in S\} \subseteq I$ is a point-finite family then there exists $S' \subseteq S$ such that $$\bigcup_{s \in S'} A_s \not \in \{B \triangle N \colon B \in \mathcal{B}, N \in I\}.$$}

		Another important result in this topic is in \cite{GL}, where the authors cite  the following result came from Prikry.
		\\\\
	\textbf{Theorem \cite{GL}}
		\textit{	Let $Y$ be a nonmeager topological space equipped with a finite regular Borel measure $\mu$ and let $\{A_s \colon s \in S\}$ be a point-finite family of subsets of $Y$ such that for any $s \in S$ $\mu(A_s) = 0$, (resp. $A_s$ is meager) but $\mu(\bigcup_{s \in S} A_s) >0$ (resp. $\bigcup_{s \in S} A_s$ is nonmeager). Then there exists $S' \subseteq S$ such that $\bigcup_{s \in S'} A_s$ is not measurable (resp. does not have the Baire property).}
	\\
		
		In \cite{PD}, the author showed that any uncountable point-finite family of analytic subsets of Polish space $X$  must contain a subfamily whose union is not analytic.

	  It seems to be a natural generalization because a partition is a special kind of point-finite family. Unfortunately, this topic cannot be extended even to so called point-countable families. 
	
	The generalization of the problem of the existence of nonmeasurable sets  for point-finite families raises further problems. Generally, it is impossible to transfer the methods used for partitions to the new situation. This requires modifying or even changing the methods used for the proving. 
	However, there is one method that turns out to be the only one that proves the existence of sets that are nonmeasurable for both partitions and point-finite families.
	
	In 2019, in \cite{FJW},there were published the results of the existence of nonmeasurable sets in the tree structures. Exact evidence is given for the structure of Marczewski and Laver. Due to the similarity of the methods used in the evidence, the existence of sets nonmeasurable for Ellentuck spaces was also obtained.
	The main method of proving here was to use Fusion Lemma which is true for tree structures and for structure in Ellentuck space. This method seems universal enough to try to generalize the results of \cite{FJW} work for point-finite families and this is the main aim of results presented in this paper.
	
	From the formal point of view the main results presented here are more general than those given in \cite{FJW}, since, as was mentioned above, a partition is a special case of a point-finite cover. From the technical point of view, the proofs in both results are also  based on Fusion Lemma, (see \cite{AT, JB, TJ1}) which cannot be omitted in proofs, but there are some nuances which are very important to emphasize. (The proofs of the results presented in Section 3 are slightely based on the proofs of results presented in \cite{FJW} but they are modified where was essential).
	Such an adaptation of the method to new assumptions being a generalization of the previous one, which is shown here, can be used in the proofs of the existence of other structures where  generalizations are possible.
	
	 It is important to emphasize that our results are not only true for the whole spaces but also for the "large" subsets.
	 
	In this paper we consider three structures: the structure of Marczewski sets (Sacks forcing), the structure of Laver structure (Laver forcing), about which we do not know if they make topology, and Ellentuck topology (Mathias forcing).
	For Marczewski structure $(s)$-sets and $ (s^0) $-sets will be equivalent to large and small sets, respectively,
	for Laver structure $(l)$-sets and $( l^0) $-sets will be equivalent to large and small sets, respectively. 
	Finally, for Ellentuck structure $CR$-sets and $ NR $ -sets are the equivalent to large and small sets, respectively.
	If we adopt the nomenclature of large and small sets for Measure and Category Theory we associate small sets with meager sets (sets of measure zero) and large sets with sets with the Baire property (sets with positive measure).
	 
	  Our results generalize the previous results for both point-finite covers and partitions, also from \cite{FJW}, which become a consequence of Theorem 1. 
	
	The paper is organized as follows. In Section 2 there are given definitions and previous results, in Section 3 there are given main results. 
	
	We use the standard terminology for the field. The definitions and facts not cited here  one can be find in \cite{RE, KK} (topology), \cite{EL} (Ellentuck topology), \cite{TJ} (set theory), \cite{JB, TJ1, JMS} (forcing).

		\section{Definitions and previous results}
	
	\subsection{Tree ideals}
	
	Let $K\subseteq \omega$ be a set, (finite or infinite).
	A set $T\subseteq K^{<\omega}$ is called a \textit{tree} iff $t\upharpoonright n \in T$ for all $t \in T$ and $n \leqslant |t|$, i.e. $T$ is closed downwards under initial segments. It is assumed that trees have no terminal nodes.
	Let $\mathbb{T}$ means a family of all tress. 
	For each $T \in \mathbb{T}$ and $t \in T$ the set 
	$$split(t, T) = |\{n \in K \colon t^\smallfrown n \in T\}|$$ denotes the number of successors of nodes in $T$.

	\begin{definition}
		A tree $T$ is called 
		\begin{enumerate}
			\item Sacks or perfect tree iff $K=\{0,1\}$ and for each $t \in T$ there is an $s \in T$ such that $t \subseteq s$ and $split(t, T) = 2$,
			\item Laver tree iff $K=\omega$ and there is $s \in T$ such that for each $t \in T$ 
				\begin{itemize}
					\item [(a)] either $t \subseteq s$ or $s \subseteq t$,	
					\item [(b)] $split(t, T)$ is infinite for each $t \in T$.  
			\end{itemize}			
		\end{enumerate}
	\end{definition}

We denote by $\mathbb{S}$ ($\mathbb{L}$) the family of all Sacks (Laver) trees, respectively.

	In the results below, we use a Laver-like structure of a tree, (which is rather close to superset or Miller trees), because the $n-$th level of nodes we divide  into $(n + 1)$ subsets of nodes and so for each $n \in \omega$.
	\\
	
	We say that $t\in Lev_n (T)$, (i.e. $t$ belongs to $n-$level of $T$) iff there are $n-$splits below $t$.
	\\
	
	Let $$[T] = \{{\color{red}s} \in K^\omega \colon \forall_{n \in \omega}\ {\color{red}s} \upharpoonright n \in T \}$$
	be the set of all infinite paths through $T$.
	\\
	Notice that $[T]$ is closed in the Baire space $K^\omega$, (see e.g. \cite{TJ}).
	
	By $stem(T)$ we mean a node ${\color{red}s} \in T$ such that $split(s, T) = 1$ and $split(t, T) >~1$ for any $s \varsubsetneq t$.
	
	The ordering on $\mathbb{S}$ is defined as follows
	$Q\leqslant T $ iff $ Q \subseteq T$
	and 
	$$Q \leqslant_n T \textrm{ iff } Q \leqslant T \textrm{ and  any node of $n$-level of $T$ is a node of $n$-level of $Q$}.$$
	\indent
	If $T \in \mathbb{L}$, then $\{{\color{red}s} \in [T] \colon stem(T) \in {\color{red}s}\}$, (i.e. the part of $T$ above the $stem(T)$),  can be enumerated as follows:
	$$s^T_0 = stem(T), s^T_1, ..., s^T_n, ...\ .$$
	Thus, we can define the ordering on $\mathbb{L}$ in the following way: let $Q, T \in \mathbb{L}$, $Q \leqslant T$ iff $Q \subseteq T$
	and
	$$Q \leqslant_n T \textrm{ iff } stem(Q) = stem(T) \textrm{ and } s^Q_i = s^T_i \textrm{ for all } i = 0, ..., n.$$

	 We say that a set $A \subseteq K^\omega$ is a \textit{$(t)-$set} iff 
		$$\forall_{T \in \mathbb{T}}\ \exists_{Q \in \mathbb{T}}\ Q \subseteq T \wedge ([Q] \subseteq A \vee [Q]\cap A =\emptyset). $$
		We say that a set $A \subseteq K^\omega$ is a \textit{$(t^0)-$set} iff 
		$$\forall_{T \in \mathbb{T}}\ \exists_{Q \in \mathbb{T}}\ Q \subseteq T \wedge  [Q]\cap A =\emptyset. $$
		
		Throughout the paper we assume that a set of trees is a $(t^0)-$set iff its set of infinite paths is a $(t^0)-$set in $K^\omega$.
		Thus, we will denote by $\mathbb{S}^0$ $(\mathbb{L}^0)$ the family of all $(s^0)-$sets ($(l^0)-$sets), respectively. 
		
		The fact that $(s^0)-$ and $(l^0)-$sets are $\sigma-$ideals in $2^\omega$ and $\omega^\omega$, respectively,  are applications of Fusion Lemma, (see Section 2.3 below).
	
For further considerations, unless otherwise stated, $\mathbb{T}$ and $\mathbb{T}^0$ mean $\sigma-$ideals of: Sacks trees, (i.e. $(s)-$ and $(s_0)-$tree, respectively) and Laver tree (i.e. $(l)-$ and $(l^0)-$tree,respectively). Then $(t), (t^0)$ and $K$ will be determined accordingly to these structures.

	\subsection{Ellentuck topology}
	
	The Ellentuck topology $[\omega]^{\omega}_{EL}$ on $[\omega]^\omega$ is generated by sets of the form
	$$[a, A] = \{B \in [A]^\omega \colon a \subset B \subseteq a \cup A\},$$
	where $a \in [\omega]^{<\omega}$ and $A \in [\omega]^\omega$. We call such sets \textit{Ellentuck sets, (shortly $EL$-sets).} Obviously $$[a, A] \subseteq [b, B] \textrm{ iff } b \subseteq a \textrm{ and } A \subseteq B.$$ 
	
	A set $M \subseteq [\omega]^\omega$ is called a \textit{completely Ramsey} set, (shortly \textit{$CR$-set}), iff for every $[a, A]$ there exists $B \in [A]^\omega$ such that $$[a, B] \subseteq M \textrm{ or } [a, B] \cap M = \emptyset.$$ 
	A set $M \subseteq [\omega]^\omega$ is called a \textit{Ramsey  null} set, (shortly \textit{$NR$-set}), if for every $[a, A]$ there exists $B \in [A]^\omega$ such that $$[a, B] \cap M = \emptyset.$$
	\\
	Notice that  it is application of Fusion Lemma (see Section 2.3) that all $NR$-sets form a $\sigma$-ideal in $[\omega]^{\omega}_{EL}$ which we denote by $\mathbb{NR}$. 
	\\

		\noindent The following fact will be used in further considerations for simplyfying the notation. 
		\\\\
		\noindent
		\textbf{Fact 1 (\cite{JB})} \textit{Let $M$ be an open and dense set (in the sense of Ellentuck topology). Then for each $A \subseteq [\omega]^\omega$ and for each $a \in [\omega]^{<\omega}$ there exists $B \subseteq [\omega]^\omega$ such that $B \subseteq A$ the set $[\emptyset, B \cup a] \subseteq M$.}

		\subsection{Fusion Lemma}
	Let $\mathbb{T}$  be the family of all trees.  
	A sequence $\{T_n\}_{n \in \omega}$ of trees such that 
	$$T_0 \geqslant_0 T_1 \geqslant_1 ... \geqslant_{n-1} T_n \geqslant_n ...$$
	is called a \textit{fusion sequence}.
	\\
	\\
	\textbf{Fact 2 (\cite{TJ1})} If $\{T_n\}_{n\in \omega}$ is a fusion sequence then $T = \bigcap_{n\in \omega}T_n$, (the fusion of $\{T_n\}_{n \in \omega}$), belongs to $\mathbb{T}$.
	\\ 

	A sequence $\{[a_n, A_n]\}_{n \in \omega}$ of $EL$-sets  is called a \textit{fusion sequence} if it is infinite and
	\\(1) $\{a_n\}_{n \in \omega}$ is a nondecreasing sequence of integers converging to infinity,
	\\(2) $A_{n+1} \in [a_n, A_n]$ for all $n \in \omega$.
	\\
	\\
	\textbf{Fact 3 (\cite{TJ1})} If $\{[a_n, A_n]\}_{n \in \omega}$ is a fusion sequence then $$[a, A] = \bigcap_{n\in \omega}[a_n, A_n] =  [{\color{red}a},\bigcap_{n\in \omega}  A_n],$$ (the fusion of $\{[a_n, A_n]\}_{n \in \omega})$, is an $EL$-set.
	
\subsection{Point-finite covers}

Let $B \in K^\omega$, where $K$ is as described in Section 2.1. 
A family $\mathcal{A} \subseteq P(B)$ is called a \textit{point-finite cover} of $B$ iff $\bigcup\mathcal{A} = B$ and $$\{A \in \mathcal{A} \colon b \in A\}$$
is finite for each $b \in B$.

Let $M \subseteq [\omega]^\omega$.
A family $\mathcal{A} \subseteq P(M)$ is called a \textit{point-finite cover} of $M$ iff $\bigcup\mathcal{A} = M$ and $$\{A \in \mathcal{A} \colon a \in A\}$$
is finite for each $a \in M$.
	
	\section{Main results}
	
	In this part we present the main results of this paper.
	We start with auxiliary lemmas which play the crucial role in Theorem 1. 
	 As was mentioned in Section 2.1,  in the construction given in Lemma 1, we rather use a Laver-like structure of a tree, because the $n-$th level of nodes we divide  into $(n + 1)$ subsets of nodes and so for each $n \in \omega$.
	
		\begin{lemma}
			\begin{enumerate}
		\item Let $A \in P(K^\omega) \setminus \mathbb{T}^0$. For any  point-finite cover $\mathcal{F}$ of $A$ consisiting of  $(t^0)$-sets and for any perfect tree $T\in \mathbb{S}$ there exists a perfect subtree $Q \leqslant T$  such that the family
		$$\mathcal{F}_{[Q]} = \{F \cap [Q] \colon F_\alpha \in \mathcal{F}\}$$ has cardinality continuum.
		\item 	Let $M \in P([\omega]^\omega)\setminus \mathbb{NR}$ be an open and dense set. For any point-finite cover $\mathcal{F}$ of $M$ consisting of $NR$-sets and for any $[a, A] \subseteq [\omega]^{\omega}_{EL} $ there exists $[b, B] \subseteq [a, A]$ such that the family $$\mathcal{F}_{[b, B]} = \{F \cap [b, B] \colon F \in \mathcal{F}\}$$ has cardinality continuum.
		\end{enumerate}
	\end{lemma}
	
	\begin{proof}
	\textbf{1.} Let $A \in P(K^\omega) \setminus \mathbb{T}^0$ be an arbitrary set and let $T \in \mathbb{T}$. Let $\mathcal{F}$ be a point-finite cover of $A$ consisting of $(t^0)-$sets. 
		We will construct inductively by $n \in \omega$ a collection of subfamilies  $\{\mathcal{F}_h \colon h\in k^n\}$, ($k=2$ for Sacks trees and $k=n$ for Laver trees,  of $\mathcal{F}$ and a collection of perfect subtrees $\{T_h \colon h \in k^n\}$ of $T$ with the following properties:
		for any distinct $h, h' \in k^n$ 
		\begin{itemize}
	
	\item[(i)] $\mathcal{F}_h \subseteq \mathcal{F}$ and $ T_h \leqslant T$;
		
	\item[(ii)] $\bigcup\{\mathcal{F}_{h} \colon h \in k^n \} = \bigcup \mathcal{F}$;
		
	\item[(iii)] $\bigcup\mathcal{F}_h \not \in \mathbb{T}^0$;
		
	\item[(iv)] $\mathcal{F}_{h} \subseteq \mathcal{F}_{g}$ and $ T_h \leqslant_n T_{g} $, i.e. $[T_h] \subseteq [T_{g}]$, for $h\cap g = g$;
		
	\item[(v)] $\mathcal{F}_{h} \cap \mathcal{F}_{h'} = \emptyset$ and $[T_h]\cap [T_{h'}] =\emptyset$;
		
	\item[(vi)] $A \cap [T_h] \subseteq \mathcal{F}_h$ and $A\cap [T_h] \not \in \mathbb{T}^0$.
	\end{itemize}

The first and the successor step are essentially the same. 	
		Assume that for some $m \in \omega$ we have constructed the families $$\{\mathcal{F}_h \colon h \in k^m \}  \textrm{ and } \{T_h \colon h \in k^m\}$$ of properties (i) - (vi).
		
		Now, fix $h \in k^m$. Denote 
		$$W_n = \{x \in A \colon |\{F \in \mathcal{F} \colon x \in F\}| =n\}.$$
		Since $\bigcup \mathcal{F}_h  \not \in \mathbb{T}^0$ then there exists $n \in \omega$ such that $$W_n \cap \bigcup \mathcal{F}_h  \not \in \mathbb{T}^0.$$
		Let $$n_h = \min\{n \in \omega \colon W_n \cap \bigcup \mathcal{F}_h \not \in \mathbb{T}_0\}.$$
		We show that for any $A_0$ and $ A_1 = \bigcup \mathcal{F}_h \setminus A_0$ such that $W_{n_h}\cap A_\varepsilon \not \in \mathbb{T}^0$ and 
		$$W_{n_h}\cap A_\varepsilon \cap (\bigcup \mathcal{F}_{h^{\smallfrown} \nu}\setminus  \bigcup \mathcal{F}_{h^{\smallfrown}(1- \nu)}) \not \in \mathbb{T}^0$$
		for all $\varepsilon, \nu \in \{0,1\}$.
		
		Consider
		$$\{\mathcal{G}\subseteq \mathcal{F}_h \colon \bigcup \mathcal{G} \cap A_0 \cap W_{n_h} \in \mathbb{T}^0\}.$$
		Then, by Ulam Theorem, (see \cite{SU} and \cite[p.86]{KK}), there are $\mathcal{G}_0$ and $ \mathcal{G}_1 = \bigcup \mathcal{F}_h \setminus \mathcal{G}_0$ and $$\bigcup \mathcal{G}_\nu \cap A_0 \cap W_{n_h} \not \in \mathbb{T}^0$$ for each $\nu \in \{0,1\}$. 
		
		If one of these intersections  belongs to $\mathbb{T}^0$ then the second one must belong to $P(K^\omega) \setminus \mathbb{T}^0$  and 
		$$A' = \bigcup \mathcal{G}_0 \cap \bigcup \mathcal{G}_1 \cap A_0 \cap W_{n_h} \not \in \mathbb{T}^0.$$ Then, we continue the division of $\mathcal{G}_0$ and $\mathcal{G}_1$ analogously as above obtaining $\mathcal{G}_{00}, \mathcal{G}_{01} = \bigcup \mathcal{F}_h \setminus \mathcal{G}_{00},\mathcal{G}_{10}, \mathcal{G}_{11} = \bigcup \mathcal{F}_h \setminus \mathcal{G}_{10}$ such that
		$$A'' = A' \cap (\bigcup \mathcal{G}_{00} \setminus \bigcup \mathcal{G}_{01}) \not \in \mathbb{T}^0$$  
		$$A''' = A' \cap (\bigcup \mathcal{G}_{01} \setminus \bigcup \mathcal{G}_{00}) \not \in \mathbb{T}^0$$
		$$A^{iv} = A'' \cap (\bigcup \mathcal{G}_{10} \setminus \bigcup \mathcal{G}_{11}) \not \in \mathbb{T}^0$$
		and
		$$A^{v} = A'' \cap (\bigcup \mathcal{G}_{11} \setminus \bigcup \mathcal{G}_{10}) \not \in \mathbb{T}^0.$$
		Now, put $\mathcal{H}_0 = \mathcal{G}_{00}\cup \mathcal{G}_{10}$ and $\mathcal{H}_1 = \mathcal{G}_{01}\cup \mathcal{G}_{11}.$
		Obviously, $\mathcal{H}_0 \cup \mathcal{H}_1 = \bigcup\mathcal{F}_h$ and
		$$A^{v} \cap (\bigcup \mathcal{G}_{11} \setminus \bigcup \mathcal{G}_{10}) \subseteq A' \cap (\bigcup \mathcal{H}_1 \setminus \bigcup \mathcal{H}_0)$$
		$$A^{v} \cap (\bigcup \mathcal{G}_{10} \setminus \bigcup \mathcal{G}_{11}) \subseteq A' \cap (\bigcup \mathcal{H}_0 \setminus \bigcup \mathcal{H}_1)$$
		$$A_0 \cap (\bigcup \mathcal{H}_1 \setminus \bigcup \mathcal{H}_0) \cap W_{n_h} \not \in \mathbb{T}^0$$
		and
		$$A_0 \cap (\bigcup \mathcal{H}_0 \setminus \bigcup \mathcal{H}_1) \cap W_{n_h} \not \in \mathbb{T}^0.$$
		If one of the above sets belongs to $\mathbb{T}^0$ then we proceed as above producing $\mathcal{H}_{00}, \mathcal{H}_{01}, \mathcal{H}_{10},\mathcal{H}_{11}$, etc. up to obtaining both sets not belonging to $\mathbb{T}^0$.
		Thus, put
		$$\mathcal{F}_{{h^{\smallfrown} 0}} = \{F \in \mathcal{F}_h \colon F \in \mathcal{H}_0\}$$
		and $$\mathcal{F}_{{h^{\smallfrown} 1}} = \mathcal{F}_h \setminus \mathcal{F}_{{h^{\smallfrown} 0}}.$$
		Obviously  $\mathcal{F}_{{h^{\smallfrown} 0}}$ and $ \mathcal{F}_{{h^{\smallfrown} 1}}$ have properties (i) - (iii).

		Now, we will construct $T_{h^{\smallfrown} 0}, T_{h^{\smallfrown} 1} \leqslant_m T_h$ of properties (iv) and (vi). 
		Let $B_h(m)$ be a set of all nodes of $m$-level of $T_h$ with the property
		$$t\in B_h(m) \textrm{ iff } \exists_{x \in [T_h]}\  (t \in x \wedge x \cap \bigcup \mathcal{F}_h \not = \emptyset). $$
		Since $split(t, T) = 2$ for any $t \in B_h(m)$ we can divide $B_h(m)$ into disjoint sets $B_{h^\smallfrown 0}(m+1), B_{h^\smallfrown 1}(m+1)$ which will be $(m+1)-$levels of trees  $T_{h^\smallfrown 0}, T_{h^\smallfrown 1} \leqslant_m T_h$ fulfilling (iv) - (vi), respectively.
		\\
		Define
		$$t\in B_{h^\smallfrown 0}(m+1) \leftrightarrow \exists_{x \in [T_h]}\  (t \in x \wedge x \cap \bigcup \mathcal{F}_{h^\smallfrown 0} \not = \emptyset). $$  
		In similar way we define  $B_{h^\smallfrown 1}(m+1)$. Notice that both $B_{h^\smallfrown 0}(m+1)$ and $B_{h^\smallfrown 1}(m+1)$ are nonempty because of (ii).
		Since now, the proof for Sacks and Laver trees runs differently.
		\\
		\textbf{For Sacks trees}. 
		Take subtrees $T_{h^\smallfrown 0}, T_{h^\smallfrown 1} \leqslant _m T_h$ such that
		$$\{x \in [T_h] \colon \exists_{t \in B_{h^\smallfrown 0}(m+1)\setminus B_{h^\smallfrown 1}(m+1)}\  t \in x\} \subseteq [T_{h^\smallfrown 0}]$$
		$$\{x \in [T_h] \colon \exists_{t \in B_{h^\smallfrown 1}(m+1)\setminus B_{h^\smallfrown 0}(m+1)}\  t \in x\} \subseteq [T_{h^\smallfrown 1}].$$
		The construction of $(m+1)$-step is complete.
		Now take $$Q = \bigcap_{n \in \omega} \bigcup_{h \in 2^n} T_h.$$
		By Fact 2, we have  $Q\in \mathbb{T}$. Then $\mathcal{F}_{[Q]}$ has the required property. 
		\\
		\textbf{For Laver trees}.
		Let $B_h(m)$ be a set of all nodes of $m$-level of $T_h$ with the property
		$$t\in B_h(m) \textrm{ iff } \exists_{x \in [T_h]}\  (t \in x \wedge x \cap \bigcup \mathcal{F}_h \not = \emptyset). $$
		Since $split(t, T) = \omega$ for any $t \in B_h(m)$ we can divide $B_h(m)$ into disjoint sets $B_{h^\smallfrown k}(m+1), k = 0,..., m$, which will be $(m+1)-$levels of Laver trees $T_{h^\smallfrown k}\leqslant_m T_h, k = 0,..., m$, fulfilling properties (iv) - (vi), respectively.
		\\
		Define
		$$t \in B_{h^\smallfrown 0}(m+1) \textrm{ iff } \exists_{x \in [T_h]}\  (t \in x \wedge x \cap \bigcup \mathcal{F}_h \not = \emptyset)$$
		and for  $t \in B_{h^\smallfrown k}(m+1)$ for $k=1,..., m$ iff
		\\
		(a) $t \not \in \{B_{h^\smallfrown l}(m+1) \colon l< k\}$ and 
		\\
		(b) $\exists_{x \in [T_h]}\  (t \in x \wedge x \cap \bigcup \mathcal{F}_h \not = \emptyset). $
		\\		
		Notice that all $B_{h^\smallfrown k}(m+1)$ are nonempty for all $k =0, ..., m$, because of (iii).
		
		Now take the  subtree $T_{h^\smallfrown k} \leqslant _m T_h, k = 0,..., m$, such that
		$$\{x \in [T_h] \colon \exists_{t \in B_{h^\smallfrown 1}(m+1)}\  t \in x\} = [T_{h^\smallfrown 1}].$$
		The construction of $(m+1)$-step is complete.
		Now, take $$Q = \bigcap_{n \in \omega} \bigcup_{h \in 2^n} T_h.$$
		By Fact 2, the tree  $Q\in \mathbb{T}$. Then $\mathcal{F}_{[Q]}$ has the required property. 
		
		\textbf{2.} Let $M \in P([\omega]^\omega) \setminus \mathbb{NR}$ be an arbitrary set. Let $\mathcal{F}$ be a point finite cover of $M$ consisting of $NR-$sets. 
		By Fact 1, we can consider $[\emptyset, B \cup a]$ instead of $[a, A]$ for some $B \subseteq A$.
		We will construct, inductively by $n \in \omega$, a collection of subfamilies $\{\mathcal{F}_h \colon h\in 2^n\}$ of $\mathcal{F}$ and a collection of $EL$-subsets $\{[a_{h^\smallfrown n, }, A_{h^\smallfrown n}] \colon h \in 2^n\}$ of $[\emptyset, B \cup a]$  with the following properties:
	for any distinct $h, h' \in 2^n$ 
	\\
	(i) $\mathcal{F}_h \subseteq \mathcal{F}$ and $[a_h, A_h] \subseteq [\emptyset, B\cup a]$;
	\\
	(ii) $\bigcup\{\mathcal{F}_{h} \colon h \in 2^n \} = \bigcup \mathcal{F}$;
	\\
	(iii) $\bigcup\mathcal{F}_h \not \in \mathbb{NR}$;
	\\
	(iv) $\mathcal{F}_{h} \subseteq \mathcal{F}_{h\upharpoonright{(n-1)}}$ and  $[a_h, A_h] \subseteq [a_{h\upharpoonright{(n-1)}}, A_{h\upharpoonright{(n-1)}}]$;
	\\
	(v) $\mathcal{F}_{h} \cap \mathcal{F}_{h'} = \emptyset$ and $a_h\cap a_{h'} =\emptyset$;
	\\
	(vi) $M \cap [a_h, A_h] \subseteq \mathcal{F}_h$ and $M\cap [a_h, A_h] \not \in \mathbb{NR}$.
	
	The first  and the successor step are essentially the same. 	
	Assume that for some $m \in \omega$ we have constructed the families $$\{\mathcal{F}_h \colon h \in 2^m \} \textrm{ and } \{[a_h, A_h] \colon h \in 2^m\}$$ of properties (i) - (vi).
	
	Now, fix $h \in 2^m$. Denote 
	$$W_n = \{x \in M \colon |\{F \in \mathcal{F} \colon x \in F\}| =n\}.$$
	Since $\bigcup \mathcal{F}_h  \not \in \mathbb{NR}$ then there exists $n \in \omega$ such that $$W_n \cap \bigcup \mathcal{F}_h  \not \in \mathbb{NR}.$$
	Let $$n_h = \min\{n \in \omega \colon W_n \cap \bigcup \mathcal{F}_h \not \in \mathbb{NR}\}.$$
	We show that for any $M_0 $ and $ M_1 = \bigcup \mathcal{F}_h \setminus M_0$ such that $W_{n_h}\cap M_\varepsilon \not \in \mathbb{NR}$ and 
	$$W_{n_h}\cap M_\varepsilon \cap (\bigcup \mathcal{F}_{h^{\smallfrown} \nu}\setminus  \bigcup \mathcal{F}_{h^{\smallfrown}(1- \nu)}) \not \in \mathbb{NR}$$
	for all $\varepsilon, \nu \in \{0,1\}$.
	
	Consider
	$$\{\mathcal{G}\subseteq \mathcal{F}_h \colon \bigcup \mathcal{G} \cap M_0 \cap W_{n_h} \in \mathbb{NR}\}.$$
	Then, by Ulam Theorem (see \cite{SU} and \cite[p.86]{KK}) there are $\mathcal{G}_0, \mathcal{G}_1 = \bigcup \mathcal{F}_h \setminus \mathcal{G}_0$ and $$\bigcup \mathcal{G}_\nu \cap M_0 \cap W_{n_h} \not \in \mathbb{NR}$$ for each $\nu \in \{0,1\}$. 
	
	If one of these intersections  belongs to $\mathbb{NR}$ then the second one must belong to $P([\omega]^\omega) \setminus \mathbb{NR}$  and 
	$$M' = \bigcup \mathcal{G}_0 \cap \bigcup \mathcal{G}_1 \cap M_0 \cap W_{n_h} \not \in \mathbb{NR}.$$ Then we continue division of $\mathcal{G}_0$ and $\mathcal{G}_1$ analogously as above obtaing $\mathcal{G}_{00}, \mathcal{G}_{01} = \bigcup \mathcal{F}_h \setminus \mathcal{G}_{00},\mathcal{G}_{10}, \mathcal{G}_{11} = \bigcup \mathcal{F}_h \setminus \mathcal{G}_{10}$ such that
	$$M'' = M' \cap (\bigcup \mathcal{G}_{00} \setminus \bigcup \mathcal{G}_{01}) \not \in \mathbb{NR}$$  
	$$M''' = M' \cap (\bigcup \mathcal{G}_{01} \setminus \bigcup \mathcal{G}_{00}) \not \in \mathbb{NR}$$
	$$M^{iv} = M'' \cap (\bigcup \mathcal{G}_{10} \setminus \bigcup \mathcal{G}_{11}) \not \in \mathbb{NR}$$
	and
	$$M^{v} = M'' \cap (\bigcup \mathcal{G}_{11} \setminus \bigcup \mathcal{G}_{10}) \not \in \mathbb{NR}.$$
	Now, put $\mathcal{H}_0 = \mathcal{G}_{00}\cup \mathcal{G}_{10}$ and $\mathcal{H}_1 = \mathcal{G}_{01}\cup \mathcal{G}_{11}.$
	Obviously $\mathcal{H}_0 \cup \mathcal{H}_1 = \bigcup\mathcal{F}_h$ and
	$$M^{v} \cap (\bigcup \mathcal{G}_{11} \setminus \bigcup \mathcal{G}_{10}) \subseteq M' \cap (\bigcup \mathcal{H}_1 \setminus \bigcup \mathcal{H}_0)$$
	$$M^{v} \cap (\bigcup \mathcal{G}_{10} \setminus \bigcup \mathcal{G}_{11}) \subseteq M' \cap (\bigcup \mathcal{H}_0 \setminus \bigcup \mathcal{H}_1)$$
	$$M_0 \cap (\bigcup \mathcal{H}_1 \setminus \bigcup \mathcal{H}_0) \cap W_{n_h} \not \in \mathbb{NR}$$
	and
	$$M_0 \cap (\bigcup \mathcal{H}_0 \setminus \bigcup \mathcal{H}_1) \cap W_{n_h} \not \in \mathbb{NR}.$$
	If one of the above sets belongs to $\mathbb{NR}$ then we proceed as above producing $\mathcal{H}_{00}, \mathcal{H}_{01}, \mathcal{H}_{10},\mathcal{H}_{11}$, etc. up to obtaining both sets not belonging to $\mathbb{NR}$.
	Thus, put
	$$\mathcal{F}_{{h^{\smallfrown} 0}} = \{F \in \mathcal{F}_h \colon F \in \mathcal{H}_0\}$$
	and $$\mathcal{F}_{{h^{\smallfrown} 1}} = \mathcal{F}_h \setminus \mathcal{F}_{{h^{\smallfrown} 0}}.$$
	Obviously  $\mathcal{F}_{{h^{\smallfrown} 0}}, \mathcal{F}_{{h^{\smallfrown} 1}}$ have properties (i) - (iii).
	
	 (In order to construct $EL$-sets $[a_{h^\smallfrown 0}, A_{h^\smallfrown 0}], [a_{h^\smallfrown 1}, A_{h^\smallfrown 1}]$ we will use construction similar to Mathias forcing~\cite{JB}).
	
	Take $[a_h, A_h]$ associated with $\mathcal{F}_h$. Enumerate all subsets of $a_h$ by $s_1, ..., s_k,$ where $ k = 2^m$. Construct the sequences of subsets of $A_h$:
	$$C^{h}_{0} \supseteq C^{h}_{1}\supseteq ... \supseteq C^{h}_{k}
	\textrm{ and }
	D^{h}_{0} \supseteq D^{h}_{1}\supseteq ... \supseteq D^{h}_{k}$$ 
	as follows: let $C^{h}_{0},  D^{h}_{0} \subseteq [A_{h}\setminus a_h]^\omega $ and $C^{h}_{0}\cap D^{h}_{0}  = \emptyset$.
	Given $C^{h}_{i}, D^{h}_{i}$ if there exist $C\subseteq C^{h}_{i}$ and $D\subseteq D^{h}_{i}$ such that
	$M\cap [s_i, C] \subseteq  \bigcup \mathcal{F}_{h\smallfrown 0}$
	and 
	$M\cap [s_i, D] \subseteq  \bigcup \mathcal{F}_{h\smallfrown 1}$
	then $C^{h}_{i+1}:=C$ and $D^{h}_{i+1}:=D$. If not, then then we take $C^{h}_{i+1}:=C^{h}_{i}$ and $D^{h}_{i+1}:=D^{h}_{i}$.
	\\
	Take
	$$[a_{h^\smallfrown 0}, A_{h^\smallfrown 0}] = [a_h \cup \{\min C^h_k\}, C_{k}^{h} \setminus \{\min C^h_k\}]$$ and 
	$$[a_{h^\smallfrown 1}, A_{h^\smallfrown 1}] = [a_h \cup \{\min D^h_k\}, D_{k}^{h} \setminus \{\min D^h_k\}].$$
	Thus, we have constructed the collection of families $\{\mathcal{F}_h \colon f \in 2^n\}$ and associated with them $[a_h, A_h] \subseteq [\omega]^{\omega}_{EL}$ fulfilling (i) - (vi) for all $n \in \omega$.
	Now take
	$$[\emptyset, B] = \bigcap_{n \in \omega} \bigcup_{h \in 2^n} [a_h, A_h].$$
	By Fact 3, the  set $[\emptyset, B]$ is an $EL$-set. Thus, we obtain that $\mathcal{F}_{[\emptyset, B]}$ has the required property.
	\end{proof}

	\begin{theorem}
		
		\begin{enumerate}
			\item 
			Let $A \in P(K^\omega)\setminus \mathbb{T}^0$ and let $\mathcal{F}$ be a point-finite cover of $A$ consisting of $(t^0)$-sets. If for each $T \in \mathbb{T}$, with $A \cap [T] \not = \emptyset$, there exists a subtree $Q \leqslant T$, (with $A \cap [Q] \not = \emptyset$),  such that $$\mathcal{F}_{[Q]} = \{F \cap [Q] \colon F \in \mathcal{F}\}$$ has cardinality continuum then 	$\bigcup \mathcal{F'}_{[Q]}$ is not a $(t)$-set for some subfamily $\mathcal{F}' \subseteq \mathcal{F}$.
			\item 	Let $M \in P([\omega]^{\omega})\setminus \mathbb{NR}$ and let $\mathcal{F}$ be a point-finite cover of $M$ consisting of $NR$-sets. If for each $EL$-set $[a, A]$, with $M \cap [a, A] \not = \emptyset$, there exists an $EL$-set $[b, B] \leqslant [a, A]$, (with $M \cap [b, B] \not = \emptyset$),  such that $$\mathcal{F}_{[b, B]} = \{F \cap [b, B] \colon F \in \mathcal{F}\}$$ has cardinality continuum then 	$\bigcup \mathcal{F'}_{[b, B]}$ is not a $CR$-set for some subfamily $\mathcal{F}' \subseteq \mathcal{F}$.
		\end{enumerate}
	\end{theorem}
	
	\begin{proof}
		For our convenience we will show the proof of the first part of theorem. The second part is similar.
		
		Enumerate 
		$$\mathbb{T}' = \{T_\alpha \in \mathbb{T} \colon A\cap [T_\alpha] \not = \emptyset, \alpha \in 2^\omega\}.$$
		Let $\mathcal{F}$ be a  point-finite cover of $A$ consisting of $(t^0)$-sets. By assumption, for each  tree $T_\alpha \in \mathbb{T}'$ there exists a subtree $Q_\alpha \leqslant T_\alpha$ such that the point-finite $$\mathcal{F}_{[Q_\alpha]} = \{F \cap [Q_\alpha] \colon F \in \mathcal{F}\}$$ has cardinality continuum.
		Hence, for each $\alpha \in 2^\omega$ we will choose distinct elements $x_\alpha, y_\alpha \in [Q_\alpha]$,  such that the families
		$$\mathcal{B}^0_\alpha = \{F \in \mathcal{F} \colon x_\alpha \in F \cap [Q_\alpha]\} \setminus \{F \in \mathcal{F} \colon F \in \{\mathcal{B}^0_\beta\cup \mathcal{B}^1_{\beta}\colon \beta < \alpha\}\}$$
			$$\mathcal{B}^1_\alpha = \{F \in \mathcal{F} \colon x_\alpha \in F \cap [Q_\alpha]\} \setminus \{F \in \mathcal{F} \colon F \in \{\mathcal{B}^0_\beta\cup \mathcal{B}^1_{\beta}\cup \mathcal{B}^0_\alpha\colon \beta < \alpha\}\}$$
		are nonempty.
		(Such choice is possible because  $\mathcal{F}_{[Q_\alpha]}$ has cardinality continuum, but each $\mathcal{F}^0_\beta$ and $\mathcal{F}^1_\beta$ are finite and $\beta<\alpha$).
		Notice that $\mathcal{B}^0_\alpha$ and $\mathcal{B}^1_\alpha$ are disjoint. 
		
		Now, let $\mathcal{B}^0 = \{\mathcal{B}^0_\alpha \colon \alpha < 2^\omega\}$ and $\mathcal{B}^1 = \{\mathcal{B}^1_\alpha \colon \alpha < 2^\omega\}$. Obviously, $\mathcal{B}^0$ and $\mathcal{B}^1$ are disjoint. 
			
		Notice that $\bigcup \mathcal{B}^{\varepsilon}$ are not $(t)$-sets for any $\varepsilon \in \{0, 1\}$.
		Indeed. Suppose that  $\bigcup \mathcal{B}^{\varepsilon}$ is a $(t)$-set for some $\varepsilon \in \{0, 1\}$. Then, there exists $Q_\alpha \leqslant T_\alpha$ such that $[Q_\alpha] \cap \bigcup \mathcal{B}^{\varepsilon} = \emptyset$. But by the construction we have that 
		$$\{F \in \mathcal{F} \colon F \cap [Q_\alpha] \cap \bigcup \mathcal{B}^{\varepsilon} \not = \emptyset\}$$ is nonempty. A contradiction.
		
		If $\bigcup \mathcal{B}^{\varepsilon}$ is not a $(t)$-set for some $\varepsilon \in \{0, 1\}$, then there exists $Q_\alpha \leqslant T_\alpha$ such that $[Q_\alpha] \subseteq \bigcup \mathcal{B}^{\varepsilon}$ and by the construction $[Q_\alpha] \cap \bigcup \mathcal{B}^{1 - \varepsilon} \not = \emptyset$ which contradicts with disjointness of families $\mathcal{B}^{0}$ and  $\mathcal{B}^{1}$. 
	\end{proof}

		\noindent
	{\sc Joanna Jureczko}
	\\
	Wroc\l{}aw University of Science and Technology,
	Faculty of Electronics, Wroc\l{}aw, Poland
	\\
	{\sl e-mail: joanna.jureczko@pwr.edu.pl}
\end{document}